\documentclass[english]{article}
\usepackage[T1]{fontenc}
\usepackage[latin9]{inputenc}
\usepackage[letterpaper]{geometry}
\usepackage{amsmath}
\usepackage{setspace}
\usepackage{amssymb}
\usepackage{esint}
\makeatletter
\usepackage{babel}

\makeatletter
\newtheorem{theorem}{Theorem}

\newtheorem{example}[theorem]{Example}

\newtheorem{proposition}[theorem]{Proposition}
\newtheorem{remark}[theorem]{Remark}

\newenvironment{proof}[1][Proof]{\textbf{#1.} }{\ \rule{0.5em}{0.5em}}
\makeatother

\makeatother

\usepackage{babel}

\begin{document}

\title{Enlargements of Filtrations and Applications}

\author{José Manuel Corcuera%
\thanks{Universitat de Barcelona, Gran Via de les Corts Catalanes, 585, E-08007
Barcelona, Spain. E-mail\texttt{: jmcorcuera@ub.edu}%
} %
\thanks{The work of J. M. Corcuera is supported by the NILS Grant. %
} \and Arturo Valdivia%
\thanks{Universitat de Barcelona, Gran Via de les Corts Catalanes, 585, E-08007
Barcelona, Spain. E-mail\texttt{: arturo.valdivia@ub.edu}%
} %
\thanks{Corresponding author.%
}}
\maketitle
\begin{abstract}
In this paper we review some old and new results about the enlargement
of filtrations problem, as well as their applications to credit risk
and insider trading problems. The enlargement of filtrations problem
consists  in the study of conditions under which a semimartingale
remains a semimartingale when the filtration is enlarged, and, in
such a case, how to find the Doob-Meyer decomposition. Filtrations
may be enlarged in different ways. In this paper we consider initial
and progressive filtration enlargements made by random variables and
processes.

\textbf{Keywords}: Credit Risk, Insider Trading, Enlargement of Filtrations

\textbf{Mathematics Subject Classification (2010)}: 60-02, 60G44,
62B10, 91-02, 91G40
\end{abstract}

\section{Enlargement of filtrations}

When considering a given filtration $\mathbb{F}=(\mathcal{F}_{t})_{t\ge0}$,
accounting for the information related to a given phenomenon, the
\emph{arrival} of new information induces the consideration of an
\emph{enlarged filtration} $\mathbb{G}=(\mathcal{G}_{t})_{t\ge0}$,
such that $\mathcal{G}_{t}\supseteq\mathcal{F}_{t}$, for each $t\ge0$.
More specifically, one considers the filtration $\mathbb{G}:=(\mathcal{G}_{t})_{t\ge0}$
defined as $\mathcal{G}_{t}:=\mathcal{F}_{t}\vee\mathcal{H}_{t}\supseteq\mathcal{F}_{t}$,
$t\ge0$, where $\mathbb{H}=(\mathcal{H}_{t})_{t\ge0}$ is assumed
to represent the new information. Traditionally, when $\mathbb{H}$
is such that $\mathcal{H}_{t}=\sigma(L)$, for $t\ge0$, and for some
random variable $L$, then it is said that $\mathbb{G}$ is an \emph{initial
enlargement }of the filtration\emph{\ }$\mathbb{F}$. Otherwise,
it is said that $\mathbb{G}$ is a \emph{progressive enlargement}
of the filtration $\mathbb{F}$\emph{.}

Some natural questions arise in a filtration enlargement setting.
For instance, \emph{In which cases an $\mathbb{F}$-semimartingale
remains a semimartingale in the enlarged filtration $\mathbb{G}$?},
and \emph{How can we compute $\mathbb{G}$-Doob-Meyer decompositions
in function of $\mathbb{F}$ and $\mathbb{H}$? }According to (1),
the beginnings of the theory of enlargement of filtrations may be
traced back to K. Itô, P.A. Mayer, and D. Williams who in the late
seventies, independently and separately from each other, posed similar
questions. So far, the study of enlargement of filtrations has been
devoted mainly to cases as $(\mathcal{H}_{t}=\sigma(\tau))_{t\ge0}$
and $(\mathcal{H}_{t}=\mathbf{1}_{\{\tau\leq t\}})_{t\ge0}$, where
$\tau$ is some stopping time. One of the drawbacks of the present
approaches is that some rather restrictive or unrealistic assumptions
has to be made on the stopping time in order to apply the approach.
For instance, in credit risk theory, $\tau$ usually represents the
\emph{default time} of some contract. In order to preclude arbitrage,
$\tau$ is assumed to be either \emph{initial time }or \emph{honest
time}%
\footnote{See section \ref{sub:credit-risk} for more details.%
}. On the other hand, only a few studies has been developed in the
general setting. For instance, it can be considered the case%
\footnote{See Section 1.2.2 in (1).%
} when $\mathbb{H}=(\mathcal{H}_{t}=$ $\sigma(J_{t}))_{t\ge0}$, for
$(J_{t}=\inf_{s\geq t}X_{s})_{t\ge0}$, being $(X_{t})_{t\ge0}$ a
3-dimensional Bessel process. Nevertheless, this case can be reduced
in fact to a case with random times, taking into account that\[
\left\{ J_{t}<a\right\} =\left\{ t<\Lambda_{a}\right\} ,\]
 where $\Lambda_{a}=\sup\{t,\; X_{t}=a\}$.

In this paper we approach the problem of initial enlargement of filtrations
by considering a general $\mathcal{F}$-measurable random variable
$L$ (as opposed to a stopping time) and looking into its law conditioned
to the -not necessarily Brownian- filtration $\mathbb{F}$. Regarding
to the first question, we give a condition under which $\mathbb{F}$-semimartingales
remain semimartingales under the enlarged filtration $\mathbb{G}$.
Regarding to the second question, we study conditions under which
it is possible to obtain $\mathbb{G}$-Doob-Meyer decompositions as
well as explicit expressions for the compensantor. For progressive
enlargement of filtrations, we present the example that motivated
this paper and which represents an enlargement done by an $\mathbb{H}$
that is not induced by a stopping time. Finally, we present applications
of the enlargement of filtration theory in the field of mathematical
finance, specifically in credit risk theory and insider trading.

\subsection{Initial enlargement of filtrations}

Consider a stochastic basis $\left(\Omega,\mathcal{F},\mathbb{F},\mathbb{P}\right)$,
where the filtration $\mathbb{F}$ is assumed to satisfy the usual
conditions\emph{. }Let $L$ be an $\mathcal{F}$-measurable random
variable with values in $(\mathbb{R},\mathcal{B}\left(\mathbb{R}\right))$.
Let $T>0$ denote a time horizon, and define $\mathcal{G}_{t}:=\cap_{T\geq s>t}\left(\mathcal{F}_{s}\vee\sigma(L)\right)$
and $\mathbb{G=}(\mathcal{G}_{t})_{t\in[0,T]}$. Notice that defining
$\mathcal{G}_{t}$ as $\cap_{T\geq s>t}\left(\mathcal{F}_{s}\vee\sigma(L)\right)$
assures the right-continuity of the filtration $\mathbb{G}$, and
therefore, that $\mathbb{G}$ satisfy the usual conditions.

\textbf{Condition A}. For all $t$, there exists a $\sigma$-finite
measure $\eta_{t}$ in $(\mathbb{R},\mathcal{B}\left(\mathbb{R}\right))$
such that $Q_{t}(\omega,\cdot)\ll\eta_{t}$ where $Q_{t}(\omega,\mathrm{d}x)$
is a regular version of $L|\mathcal{F}_{t}.$

Notice that \textbf{Condition A }is satisfied in the case that $L$
takes a countable number of values and the case when $L$ is independent
of $\mathcal{F}_{\infty}$, just taking $\eta_{t}$ the law of $L.$
An example where \textbf{Condition A} is not satisfied is the following
example:

\begin{example} \label{First}Let $L$ be the $n$-th jump of a Poisson
process $(N_{t})_{t\in[0,T]}$ with intensity $\lambda$ and $\mathbb{F}$
the filtration generated by $(N_{t})_{t\in[0,T]}$, then\[
\mathbb{P}\{L>x|\mathcal{F}_{t}\}=\mathbf{1}_{\{N_{x}<n,N_{t}\geq n\}}+\mathbf{1}_{\{N_{t}<n\}}\int_{(x-t)_{+}}^{\infty}\frac{\lambda e^{-\lambda u}\left(\lambda u\right)^{n-N_{t}-1}}{(n-N_{t}-1)!}\mathrm{d}u,\]
 then the conditional probability cannot be dominated by a non random
measure. \end{example}

\begin{theorem} Under \textbf{Condition A} any $X,$ $\mathbb{F}$-semimartingale
is a $\mathbb{G}$-semimartingale.

\begin{proof} The proof is based in the characterization of semimartingales
of Bichteler-Dellacherie-Mokobodzki: Let $X$ be a càdlàg process,
adapted to a filtration $\mathbb{H}$ define the class of predictable
processes\[
\xi_{\rho,t}:=\left\{ f=\sum_{i=1}^{n-1}f_{i}\mathbf{1}_{(s_{i},s_{i+1}]},0=s_{0}<...s_{n}\leq t,\text{ }f_{i}\in H_{s_{i}-},\ |f_{i}|<\rho\right\} \]
 define \[
\int_{0}^{t}f_{s}\mathrm{d}X_{s}:=\sum_{i=1}^{n-1}f_{i}(X_{s_{i+1}}-X_{s_{i}})\]
 and, for $Z\in L^{1}$ \[
\alpha_{\rho,t}^{X}(Z,\mathbb{H})=\sup_{f\in\xi_{\rho,t}}\mathbb{E}\left[|Z|\left(1\wedge\int_{0}^{t}f_{s}dX_{s}\right)\right],\]
 then, \begin{eqnarray*}
X & \in & S(\mathbb{H})\Longleftrightarrow\lim_{\rho\rightarrow0}\alpha_{\rho,t}^{X}(1,\mathbb{H})=0\\
X & \in & S(\mathbb{H})\Longrightarrow\lim_{\rho\rightarrow0}\alpha_{\rho,t}^{X}(Z,\mathbb{H})=0.\end{eqnarray*}
 This result together with the fact that any $f,$ $\mathcal{G}_{i}$-measurable
function, can be written as $f=g(\omega)^{L(\omega)},$ where $g(\omega)^{x}$
is $\mathcal{B}(\mathbb{R})\otimes\mathcal{F}_{i}$-measurable, allows
to get the result. \end{proof} \end{theorem}

\begin{proposition} \textbf{Condition A} is equivalent to $Q_{t}(\omega,\mathrm{d}x)\ll\eta(\mathrm{d}x)$
where $\eta$ is the law of $L.$ \end{proposition}

\begin{proof} By \textbf{Condition A} we have that $Q_{t}(\omega,\mathrm{d}x)=q_{t}^{x}(\omega)\eta_{t}(\mathrm{d}x),$
where $q_{t}^{x}(\omega)$ is $\mathcal{B}(\mathbb{R})\otimes\mathcal{F}_{t}$-measurable
then we can write $Q_{t}(\omega,\mathrm{d}x)=\hat{q}_{t}^{x}(\omega)\eta(\mathrm{d}x)$
with $\hat{q}_{t}^{x}(\omega)=\frac{q_{t}^{x}(\omega)}{\mathbb{E}(q_{t}^{x}(\omega))}.$
\end{proof}

\begin{proposition} Under \textbf{Condition A} there exists $q_{t}^{x}(\omega)$
$\mathcal{B}(\mathbb{R})\otimes\mathcal{F}_{t}$-measurable such that\begin{equation}
Q_{t}(\omega,\mathrm{d}x)=q_{t}^{x}(\omega)\eta(\mathrm{d}x)\label{eq: Jacod_lemma_1.8}\end{equation}

and, for fixed $x$, $(q_{t}^{x})_{t\in[0,T]}$ is an $\mathbb{F}$-martingale.
\end{proposition}

\begin{proof} See (2) Lemma 1.8. \end{proof}

\begin{theorem} Let $(M_{t})_{t\in[0,T]}$ be a continuous local
$\mathbb{F}$-martingale and $k_{t}^{x}(\omega)$ such that \[
\langle q^{x},M\rangle_{t}=\int_{0}^{t}k_{s}^{x}q_{s-}^{x}\mathrm{d}\langle M,M\rangle_{s}\]
 then \[
M_{\cdot}-\int_{0}^{\cdot}k_{s}^{L}\mathrm{d}\langle M,M\rangle_{s},\]
 is a $\mathbb{G}$-martingale. \end{theorem}

\begin{proof} Except for a localization procedure (see details in
(2) Theorem 2.1) the proof is the following: let $s<t$, $Z\in\mathcal{F}_{s}$,
and $g$ be abounded Borel function, then\begin{eqnarray*}
\mathbb{E}\left[Zg(L)(M_{t}-M_{s})\right] & = & \mathbb{E}\left[\mathbb{E}\left[Zg(L)(M_{t}-M_{s})|\mathcal{F}_{t}\right]\right]\\
 & = & \mathbb{E}\left[Z(M_{t}-M_{s})\mathbb{E}\left[g(L)|\mathcal{F}_{t}\right]\right]\\
 & = & \int_{\mathbb{R}}g(x)\eta(\mathrm{d}x)\mathbb{E}\left[Z(M_{t}-M_{s})q_{t}^{x}\right]\\
 & = & \int_{\mathbb{R}}g(x)\eta(\mathrm{d}x)\mathbb{E}\left[Z(M_{t}q_{t}^{x}-M_{s}q_{s}^{x})\right]\\
 & = & \int_{\mathbb{R}}g(x)\eta(\mathrm{d}x)\mathbb{E}\left[Z(\langle M,q^{x}\rangle_{t}-\langle M,q^{x}\rangle_{s})\right]\\
 & = & \int_{\mathbb{R}}g(x)\eta(\mathrm{d}x)\mathbb{E}\left[Z\left(\int_{s}^{t}k_{u}^{x}q_{u-}^{x}\mathrm{d}\langle M,M\rangle_{u}\right)\right]\\
 & = & \mathbb{E}\left[Zg(L)\left(\int_{s}^{t}k_{u}^{x}q_{u-}^{x}\mathrm{d}\langle M,M\rangle_{u}\right)\right]\end{eqnarray*}
 \end{proof}

\begin{example} \label{Back}Let $T=1$, and take $(M_{t}:=B_{t},\; t\in[0,T])$
where $(B_{t})_{t\in[0,T]}$ is a standard Brownian motion, and take
$L:=B_{T}$. It follows easily from \eqref{eq: Jacod_lemma_1.8} that\[
q_{t}^{x}(\omega)\thicksim\frac{1}{(T-t)^{1/2}}\exp\left(-\frac{1}{2(T-t)}(B_{t}(\omega)-x)^{2}+\frac{x^{2}}{2}\right).\]
 Applying Itô's formula we get \[
\mathrm{d}q_{t}^{x}=q_{t}^{x}\frac{x-B_{t}}{T-t}\mathrm{d}B_{t},\]
 hence $k_{t}^{x}=\frac{x-B_{t}}{T-t}$ and \[
B_{\cdot}-\int_{0}^{\cdot}\frac{B_{T}-B_{s}}{T-s}\mathrm{d}s\]
 is a $\mathbb{G}:=\mathbb{F}^{B}\vee\sigma(B_{T})$ martingale. Note
that, by the Lévy theorem, $B_{\cdot}-\int_{0}^{\cdot}\frac{B_{T}-B_{s}}{T-s}\mathrm{d}s$
is a (standard) $\mathbb{G}$-Brownian motion and, since $B_{T}$
is $\mathcal{G}_{0}$-measurable, it is independent of $(W_{t})_{t\in[0,T]}$.
\end{example}

\begin{example} \label{G2}Note that if the filtration $\mathbb{F}$
is that generated by a Brownian motion, $(B_{t})_{t\in[0,T]},$ then
for any $\mathbb{F}$-martingale $(M_{t})_{t\in[0,T]}$ we have \[
\mathrm{d}M_{t}=\sigma_{t}\mathrm{d}B_{t}\]
 and so \[
\mathrm{d}\langle M,M\rangle_{t}=\sigma_{t}^{2}\mathrm{d}t.\]
 Also, assuming that \[
q_{t}^{x}\left(\omega\right)=h_{t}^{x}(B_{t})\]
 and $h\in C^{1,2}$ we will have that \[
\mathrm{d}q_{t}^{x}=\partial h_{t}^{x}(B_{t})\mathrm{d}B_{t},\]
 and \[
k_{t}^{x}=\frac{\partial\log h_{t}^{x}(B_{t})}{\sigma_{t}}.\]
 \end{example}

\begin{example} \label{G3}The previous example can be generalized
as follows: let $(Y_{t})_{t\in[0,T]}$ be the following Brownian semimartingale
\[
Y_{t}=Y_{0}+\int_{0}^{t}\sigma(Y_{s})\mathrm{d}B_{s}+\int_{0}^{t}b(Y_{s})\mathrm{d}s,\]
 and assume that \[
Y_{T}|\mathcal{F}_{t}\thicksim\pi(T-t,Y_{t},x)\mathrm{d}x,\]
 with $\pi$ smooth. We know that $\left(\pi(T-t,Y_{t},x)\right)_{t}$
is an $\mathbb{F}$-martingale, then \[
\mathrm{d}\pi(T-t,Y_{t},x)=\frac{\partial\pi}{\partial y}(T-t,Y_{t},x)\sigma(Y_{t})\mathrm{d}B_{t}\]
 and by the Jacod theorem\[
\int_{0}^{\cdot}\sigma(Y_{s})\mathrm{d}B_{s}-\int_{0}^{\cdot}\frac{\partial\log\pi}{\partial y}(T-s,Y_{s},Y_{T})\sigma^{2}(Y_{s})\mathrm{d}s\]
 is an $\mathbb{F\vee\sigma}(Y_{T})$-martingale, and we can write
\[
Y_{\cdot}=Y_{0}+\int_{0}^{\cdot}\sigma(Y_{s})\mathrm{d}\tilde{B}_{s}+\int_{0}^{\cdot}b(Y_{s})\mathrm{d}s+\int_{0}^{\cdot}\frac{\partial\log\pi}{\partial y}(T-s,Y_{s},Y_{T})\sigma^{2}(Y_{s})\mathrm{d}s,\]
 where $(\tilde{B}_{t})_{t\in[0,T]}$ is an $\mathbb{F\vee\sigma}(Y_{T})$-Brownian
motion. \end{example}

\bigskip{}
 We have also a similar result for locally bounded martingales.

\begin{theorem} Let $(M_{t})_{t\in[0,T]}$ be an $\mathbb{F}$-local
martingale locally bounded. Then, there exist $k_{t}^{x}(\omega)$
such that \[
\langle q^{x},M\rangle_{\cdot}=\int_{0}^{\cdot}k_{s}^{x}q_{s-}^{x}\mathrm{d}\langle M,M\rangle_{s},\]
 and \[
M_{\cdot}-\int_{0}^{\cdot}k_{s}^{L}\mathrm{d}\langle M,M\rangle_{s}\]
 is a $\mathbb{G}$-martingale. \end{theorem}

Suppose that any $\mathbb{F}$-local martingale admits a representation
of the form \[
M_{t}=M_{0}+\sum_{\left(n\right)}\int_{0}^{t}K_{s}^{n}\mathrm{d}X_{s}^{n}+\int_{0}^{t}\int_{\mathbb{R}}W(\omega,x,s)\left(Q(\omega,\mathrm{d}x,\mathrm{d}s)-\nu(\omega,\mathrm{d}x,\mathrm{d}s)\right)\]
 where $\left(X^{n}\right)$ are continuous local martingales pairwise
orthogonal, assume that $q^{x}$ admits the representation \[
q_{t}^{x}=q_{0}^{x}+\sum_{\left(n\right)}\int_{0}^{t}k_{s}^{n,x}q_{s-}^{x}\mathrm{d}X_{s}^{n}+\int_{0}^{t}\int_{\mathbb{R}}q_{s-}^{x}U_{s}^{x}\left(Q(\cdot,\mathrm{d}x,\mathrm{d}s)-\nu(\cdot,\mathrm{d}x,\mathrm{d}s)\right),\]
 then \[
M_{t}-\sum_{\left(n\right)}\int_{0}^{t}K_{s}^{n}k_{s}^{n,L}\mathrm{d}\langle X^{n},X^{n}\rangle_{s}-\int_{0}^{t}\int_{\mathbb{R}}W(\cdot,x,s)U_{s}^{L}\nu(\cdot,\mathrm{d}x,\mathrm{d}s)\]
 is a $\mathbb{G}$-martingale with continuous part,\[
M_{0}+\sum_{\left(n\right)}\int_{0}^{t}K_{s}^{n}\mathrm{d}X_{s}^{n}-\sum_{\left(n\right)}\int_{0}^{t}K_{s}^{n}k_{s}^{n,L}\mathrm{d}\langle X^{n},X^{n}\rangle_{s}\]
 and jump part\[
\int_{0}^{t}\int_{\mathbb{R}}W(\omega,x,s)\left(Q(\omega,\mathrm{d}x,\mathrm{d}s)-(1+U_{s}^{L})\nu(\omega,\mathrm{d}x,\mathrm{d}s)\right).\]

\begin{example} \label{Poisson}Consider the Poisson process $(N_{t})_{t\in[0,1]}$
of intensity $\lambda$, as well as the filtration $\mathbb{F}:=\mathbb{F}^{N}$
generated by it. Let \[
M_{t}=N_{t}-\lambda t,\quad t\in[0,1],\]
 and $L=N_{1}$, then \begin{eqnarray*}
Q_{t}(\cdot,\mathrm{d}k) & = & \mathbb{P}\{N_{1}=k|\mathcal{F}_{t}\}=\mathbb{P}\{N_{1-t}=k-N_{t}\}\\
 & = & e^{-\lambda(1-t)}\frac{\left(\lambda(1-t)\right)^{k-N_{t}}}{\left(k-N_{t}\right)!},\end{eqnarray*}
 and \begin{eqnarray*}
q_{t}^{k} & = & \frac{e^{-\lambda(1-t)}\frac{\left(\lambda(1-t)\right)^{k-N_{t}}}{\left(k-N_{t}\right)!}}{e^{-\lambda}\frac{\lambda^{k}}{k!}}\\
 & = & \frac{e^{\lambda t}\lambda^{-N_{t}}(1-t)^{k-N_{t}}k!}{\left(k-N_{t}\right)!}.\end{eqnarray*}
 Now, if there is a jump at $t$: \[
U_{t}^{k}=\frac{q_{t}^{k}-q_{t-}^{k}}{q_{t-}^{k}}=\frac{k-N_{t-}}{\lambda(1-t)}-1,\]
 so \[
1+U_{t}^{L}=\frac{N_{1}-N_{t-}}{\lambda(1-t)},\]
 and \[
N_{t}-\int_{0}^{t}\frac{N_{1}-N_{s}}{1-s}\mathrm{d}s,0\leq t<1\]
 is a $\mathbb{G}$-martingale. \end{example}

\begin{remark} A more general result, concerning Lévy processes,
can be obtained by using the characteristic function instead of the
conditional density. Let $(Z_{t})_{t\ge0}$ be a Lévy process with
characteristic function \[
\mathbb{E}[e^{i\theta Z_{t}}]=e^{t\psi(\theta)}.\]
 Let $0\leq s\leq u\leq t\leq T$. Then, in virtue of the independence
of the increments of $(Z_{t})_{t\ge0}$, we have that the following
chain of equations \begin{eqnarray*}
\mathbb{E}[e^{i\theta Z_{T}}(Z_{t}-Z_{u})h_{s}] & = & \mathbb{E}[e^{i\theta(Z_{T}-Z_{t})}]\mathbb{E}[e^{i\theta(Z_{t}-Z_{u})}(Z_{t}-Z_{u})]\mathbb{E}[e^{i\theta Z_{u}}h_{s}]\\
 & = & \mathbb{E}[e^{i\theta(Z_{T}-Z_{t})}]\left(\frac{1}{i}\mathbb{E}[e^{i\theta(Z_{t}-Z_{u})}]\partial_{\theta}\log\mathbb{E}[e^{i\theta(Z_{t}-Z_{u})}]\right)\mathbb{E}[e^{i\theta Z_{u}}h_{s}]\\
 & = & \frac{1}{i}\mathbb{E}[e^{i\theta Z_{T}}h_{s}]\partial_{\theta}\log\mathbb{E}[e^{i\theta(Z_{t}-Z_{u})}]\\
 & = & \frac{1}{i}(t-u)\mathbb{E}[e^{i\theta Z_{T}}h_{s}]\partial_{\theta}\psi(\theta).\end{eqnarray*}
 And, consequently, \[
\mathbb{E}\left[e^{i\theta Z_{T}}h_{s}\frac{Z_{t}-Z_{u}}{t-u}\right]=\frac{1}{i}\mathbb{E}[e^{i\theta Z_{T}}h_{s}]\partial_{\theta}\psi(\theta).\]
 Since the right hand side of previous equation does not depend on
$t$, we have \[
\mathbb{E}\left[e^{i\theta Z_{T}}h_{s}\frac{Z_{T}-Z_{u}}{T-u}\right]=\frac{1}{i}\mathbb{E}[e^{i\theta Z_{T}}h_{s}]\partial_{\theta}\psi(\theta).\]
 Thus, by integrating with respect $u$ and applying Fubini we obtain
\begin{eqnarray*}
\mathbb{E}\left[e^{i\theta Z_{T}}h_{s}\int_{s}^{t}\frac{Z_{T}-Z_{u}}{T-u}\mathrm{d}u\right] & = & \frac{1}{i}(t-s)\mathbb{E}[e^{i\theta Z_{T}}h_{s}]\partial_{\theta}\psi(\theta)\\
 & = & \mathbb{E}[e^{i\theta Z_{T}}h_{s}(Z_{t}-Z_{s})].\end{eqnarray*}
 In certain non-homogeneous cases we can use a similar argument. Assume
that \[
\mathbb{E}[e^{i\theta(Z_{t}-Z_{u})}]=e^{g(t,u)\psi(\theta)},\]
 then \[
\mathbb{E}\left[e^{i\theta Z_{T}}h_{s}\int_{s}^{t}\frac{Z_{T}-Z_{u}}{g(T,u)}(-1)\partial_{u}g(t,u)\mathrm{d}u\right]=\mathbb{E}[e^{i\theta Z_{T}}h_{s}(Z_{t}-Z_{s})],\]
 provided that $\int_{s}^{t}\frac{Z_{T}-Z_{u}}{g(T,u)}(-1)\partial_{u}g(t,u)\mathrm{d}u$
is well defined. \end{remark}

\textbf{Condition A} is a sufficient condition that allow us to find
the Doob decomposition in the enlarged filtration. However we have
the following proposition that allows us to obtain the Doob-Meyer
decomposition without requiring \textbf{Condition A}:

\begin{proposition} \label{PropNotJac}With the notations above,
assume that there exist $\alpha_{t}^{x}(\omega)$ such that \[
\langle\int_{a}^{\infty}Q_{t}(\cdot,\mathrm{d}x),M\rangle=\int_{0}^{t}\int_{a}^{\infty}\alpha_{s}^{x}Q_{s-}(\cdot,\mathrm{d}x)\mathrm{d}\langle M,M\rangle_{s},\text{ for all }a\in\mathbb{R},\]
 then \[
M_{\cdot}-\int_{0}^{\cdot}\alpha_{s}^{x}\mathrm{d}\langle M,M\rangle_{s},\]
 is a $\mathbb{G}$-martingale. \end{proposition}

\begin{proof} For every $Z\in\mathcal{F}_{s}$we have that \begin{eqnarray*}
\mathbb{E}[Z\mathbf{1}_{\{L>a\}}(M_{t}-M_{s})] & = & \mathbb{E}\left[\mathbb{E}[Z\mathbf{1}_{\{L>a\}}(M_{t}-M_{s})|\mathcal{F}_{t}]\right]\\
 & = & \mathbb{E}\left[Z(M_{t}-M_{s})\mathbb{E}[\mathbf{1}_{\{L>a\}}|\mathcal{F}_{t}]\right]\\
 & = & \mathbb{E}\left[Z(M_{t}-M_{s})\int_{a}^{\infty}Q_{t}(\omega,\mathrm{d}x)\right]\\
 & = & \mathbb{E}\left[Z\left(M_{t}\int_{a}^{\infty}Q_{t}(\omega,\mathrm{d}x)-M_{s}\int_{a}^{\infty}Q_{s}(\omega,\mathrm{d}x)\right)\right]\\
 & = & \mathbb{E}\left[Z\left(\langle M,\int_{a}^{\infty}Q_{\cdot}(\omega,\mathrm{d}x)\rangle_{t}-\langle M,\int_{a}^{\infty}Q_{\cdot}(\omega,\mathrm{d}x)\rangle_{s}\right)\right]\\
 & = & \mathbb{E}\left[Z\int_{s}^{t}\int_{a}^{\infty}\alpha_{u}^{x}Q_{u-}(\omega,\mathrm{d}x)\mathrm{d}\langle M,M\rangle_{u}\right]\\
 & = & \mathbb{E}\left[Z\mathbf{1}_{\{L>a\}}\int_{s}^{t}\alpha_{u}^{L}\mathrm{d}\langle M,M\rangle_{u}\right]\end{eqnarray*}
 \end{proof}

\begin{example} \label{G4}Let $(B_{t})_{t\ge0}$ be a Brownian motion
and take $\tau=\inf\{t>0,B_{t}=-1\}$. It is well known that\[
\mathbb{P}\{\tau\leq s|\mathcal{F}_{t}\}=2\Phi\left(-\frac{1+B_{t}}{\sqrt{s-t}}\right)\mathbf{1}_{\{\tau\wedge s>t\}}+\mathbf{1}_{\{s<\tau\wedge t\}},\]
 where $\Phi$ is the c.d.f. of a standard normal distribution. Then
in $t<s\wedge\tau$ we have, by Itô's formula, \[
\mathbb{P}\{\tau\leq s|\mathcal{F}_{t}\}=2\Phi\left(-\frac{1}{\sqrt{s}}\right)+\sqrt{\frac{2}{\pi}}\int_{0}^{t}\frac{1}{\sqrt{s-u}}e^{-\frac{\left(1+B_{u}\right)^{2}}{2\left(s-u\right)}}\mathrm{d}B_{u},\]
 so\[
\mathrm{d}\langle\mathbb{P}\{\tau\leq s|\mathcal{F}_{\cdot}\},B\rangle_{t}=\sqrt{\frac{2}{\pi}}\frac{1}{\sqrt{s-t}}e^{-\frac{\left(1+B_{t}\right)^{2}}{2\left(s-t\right)}}\mathrm{d}t,\]
 and \begin{eqnarray*}
 &  & \alpha_{t}^{s}Q_{t}(\cdot,\mathrm{d}s)\\
 & = & \frac{\partial}{\partial s}\left(\sqrt{\frac{2}{\pi}}\frac{1}{\sqrt{s-t}}e^{-\frac{\left(1+B_{t}\right)^{2}}{2\left(s-t\right)}}\right)=\frac{1}{\sqrt{2\pi}}\left(\frac{1}{\sqrt{\left(s-t\right)^{3}}}-\frac{\left(1+B_{t}\right)^{2}}{\sqrt{\left(s-t\right)^{5}}}\right)e^{-\frac{\left(1+B_{t}\right)^{2}}{2\left(s-t\right)}},\end{eqnarray*}
 finally \[
Q_{t}(\cdot,\mathrm{d}s)=\frac{\partial}{\partial s}\mathbb{P}\{\tau\leq s|\mathcal{F}_{t}\}=\frac{e^{-\frac{\left(1+B_{t}\right)^{2}}{2\left(s-t\right)}}}{\sqrt{2\pi}\sqrt{\left(s-t\right)^{3}}}\left(1+B_{t}\right),\]
 and \[
\alpha_{t}^{s}=\frac{\frac{\partial}{\partial s}\left(\sqrt{\frac{2}{\pi}}\frac{1}{\sqrt{s-t}}e^{-\frac{\left(1+B_{t}\right)^{2}}{2\left(s-t\right)}}\right)}{\frac{\partial}{\partial s}\mathbb{P}\{\tau\leq s|\mathcal{F}_{t}\}}=\frac{1}{1+B_{t}}-\frac{1+B_{t}}{s-t}.\]
 Consequently \[
B_{t}-\int_{0}^{t\wedge\tau}\left(\frac{1}{1+B_{s}}-\frac{1+B_{s}}{\tau-s}\right)\mathrm{d}s,\quad t\geq0,\]
 is a $\mathbb{G}$-martingale. \end{example}

\begin{proposition} If we have a random measure $P_{t}^{(1)}(\omega,\mathrm{d}x)$
and a finite deterministic measure $m(\mathrm{d}t)$ such that \[
\mathbb{E}[Zg(L)(M_{t}-M_{s})]=\mathbb{E}\left[\int_{s}^{t}\int_{\mathbb{R}}Zg(x)P_{t}^{(1)}(\omega,\mathrm{d}x)m(\mathrm{d}t)\right]\]
 and $P_{t}^{(1)}(\omega,\mathrm{d}x)=\alpha_{t}^{x}(\omega)Q_{t}(\omega,\mathrm{d}x)$
then \[
M_{\cdot}-\int_{0}^{\cdot}\alpha_{s}^{L}m(\mathrm{d}s)\]
 is a $\mathbb{G}$-martingale. \end{proposition}

\begin{proof} \begin{eqnarray*}
\mathbb{E}[Zg(L)(M_{t}-M_{s})] & = & \mathbb{E}\left[\int_{s}^{t}\int_{\mathbb{R}}Zg(x)P_{u}^{(1)}(\omega,\mathrm{d}x)m(\mathrm{d}u)\right]\\
 & = & \mathbb{E}\left[\int_{s}^{t}\int_{\mathbb{R}}Zg(x)\alpha_{u}^{x}(\omega)Q_{u}(\omega,\mathrm{d}x)m(\mathrm{d}u)\right]\\
 & = & \mathbb{E}\left[\int_{s}^{t}Z\mathbb{E}[g(L)\alpha_{u}^{L}|\mathcal{F}_{t}]m(\mathrm{d}u)\right]\\
 & = & \mathbb{E}\left[\int_{s}^{t}Zg(L)\alpha_{u}^{L}m(\mathrm{d}u)\right]\end{eqnarray*}
 \end{proof}

This expression is more appropriate to treat cases like Example \ref{First}.
Suppose that $L$ is a random time and that there exist $P_{u}^{(1)}(\omega,\mathrm{d}x)$
and $P_{u}^{(2)}(\omega,\mathrm{d}x)$ such that\begin{eqnarray*}
\mathbb{E}[Zg(L)\mathbf{1}_{\{L<s\}}(M_{t}-M_{s})] & = & \mathbb{E}\left[\int_{s}^{t}\int_{0}^{s}Zg(x)P_{u}^{(1)}(\omega,\mathrm{d}x)m_{1}(\mathrm{d}u)\right]\\
\mathbb{E}[Zg(L)\mathbf{1}_{\{L>t\}}(M_{t}-M_{s})] & = & \mathbb{E}\left[\int_{s}^{t}\int_{t}^{\infty}Zg(x)P_{u}^{(2)}(\omega,\mathrm{d}x)m_{2}(\mathrm{d}u)\right],\end{eqnarray*}
 and that $P_{u}^{(i)}(\omega,\mathrm{d}x)=\alpha_{u}^{x(i)}(\omega)Q_{t}(\omega,\mathrm{d}x),$
$i=1,2,$ then it is easy to see, by decomposing $M_{t}-M_{s}$ as
sum of increments, that \[
M_{t}-\int_{0}^{t}\mathbf{1}_{\{L<u\}}\alpha_{u}^{L(1)}m_{1}(\mathrm{d}u)-\int_{0}^{t}\mathbf{1}_{\{L>u\}}\alpha_{u}^{L(2)}m_{2}(\mathrm{d}u)-\Delta M_{L}\mathbf{1}_{\{L\leq t\}},\; t\in[0,T],\]
 is a $\mathbb{G}$-martingale.

\begin{example} Consider Example \ref{First}. $L\equiv T_{n},$\begin{eqnarray*}
\mathbb{E}[Zg(T_{n})\mathbf{1}_{\{T_{n}<s\}}(N_{t}-N_{s})] & = & \mathbb{E}[(N_{t}-N_{s})]\mathbb{E}[Z\mathbb{E}[g(T_{n})\mathbf{1}_{\{T_{n}<s\}}]]\\
 & = & \mathbb{E}\left[\int_{s}^{t}\int_{0}^{s}\lambda Zg(x)Q_{t}(\omega,\mathrm{d}x)\mathrm{d}u\right]\end{eqnarray*}
 so, \[
P_{t}^{(1)}(\omega,\mathrm{d}x)=\lambda Q_{t}(\omega,\mathrm{d}x).\]
 and\begin{eqnarray*}
 &  & \mathbb{E}[Zg(T_{n})\mathbf{1}_{\{T_{n}>t\}}(N_{t}-N_{s})]\\
 & = & \mathbb{E}\left[(N_{t}-N_{s})Z\mathbb{E}[g(T_{n})\mathbf{1}_{\{T_{n}>t\}}|\mathcal{F}_{t}]\right]\\
 & = & \mathbb{E}\left[(N_{t}-N_{s})Z\int_{t}^{\infty}\frac{\lambda e^{-\lambda(x-t)}\left(\lambda(x-t)\right)^{n-N_{t}-1}}{(n-N_{t}-1)!}g(x)\mathrm{d}x\right]\\
 & = & \mathbb{E}\left[(N_{t}-N_{s})Z\int_{t}^{\infty}\frac{\lambda e^{-\lambda(x-t)}\left(\lambda(x-t)\right)^{n-\left(N_{t}-N_{s}\right)-1-N_{s}}}{(n-\left(N_{t}-N_{s}\right)-1-N_{s})!}g(x)\mathrm{d}x\right]\\
 & = & \mathbb{E}\left[Z\sum_{k=1}^{n-1-N_{s}}\int_{t}^{\infty}\frac{\lambda e^{-\lambda(x-t)}\left(\lambda(x-t)\right)^{n-k-1-N_{s}}}{(n-k-1-N_{s})!}\frac{\left(\lambda(t-s)\right)^{k}e^{-\lambda(t-s)}}{\left(k-1\right)!}g(x)\mathrm{d}x\right]\\
 & = & \lambda(t-s)\mathbb{E}\left[Z\sum_{k=0}^{n-2-N_{s}}\int_{t}^{\infty}\frac{\lambda e^{-\lambda(x-t)}\left(\lambda(x-t)\right)^{n-k-2-N_{s}}}{(n-k-2-N_{s})!}\frac{\left(\lambda(t-s)\right)^{k}e^{-\lambda(t-s)}}{k!}g(x)\mathrm{d}x\right]\\
 & = & (t-s)\mathbb{E}\left[Z\mathbb{E}\left[\left.\frac{n-N_{t}-1}{T_{n}-t}g(T_{n})\mathbf{1}_{\{T_{n}>t\}}\right|\mathcal{F}_{t}\right]\right]\\
 & = & \mathbb{E}\left[Z\int_{s}^{t}\int_{t}^{\infty}\frac{n-N_{t}-1}{x-t}g(x)Q_{t}(\omega,\mathrm{d}x)\mathrm{d}t\right],\end{eqnarray*}
 therefore\[
P_{t}^{(2)}(\omega,\mathrm{d}x)=\frac{n-N_{t}-1}{x-t}Q_{t}(\omega,\mathrm{d}x).\]
 Consequently\[
N_{t}-\lambda(t-T_{n}\wedge t)-\int_{0}^{T_{n}\wedge t}\frac{n-N_{u}-1}{T_{n}-u}\mathrm{d}u-\mathbf{1}_{\{T_{n}\leq t\}},\; t\ge0,\]
 is a $\mathbb{G}$-martingale. We also can write \[
N_{t}-\lambda(t-T_{n}\wedge t)-\int_{0}^{T_{n}\wedge t}\frac{N_{T_{n}-}-N_{u}}{T_{n}-u}\mathrm{d}u-\mathbf{1}_{\{T_{n}\leq t\}},\]
 we can compare to Example \ref{Poisson}, note also that there is
not jump at $T_{n}$. We can use Proposition \ref{PropNotJac} instead.
If we take $t<T_{n}\wedge s$\begin{eqnarray*}
\mathbb{P}\{T_{n}>s|\mathcal{F}_{t}\} & = & \int_{(s-t)_{+}}^{\infty}\frac{\lambda e^{-\lambda u}\left(\lambda u\right)^{n-N_{t}-1}}{(n-N_{t}-1)!}\mathrm{d}u\\
 & = & \int_{s}^{\infty}\frac{\lambda e^{-\lambda(u-t)}\left(\lambda(u-t)\right)^{n-N_{t}-1}}{(n-N_{t}-1)!}\mathrm{d}u=\int_{s}^{\infty}Q_{t}(\cdot,\mathrm{d}u),\end{eqnarray*}
 then, with $(M_{t}:=N_{t}-\lambda t,\; t\in[0,T])$ we have\[
\langle\int_{s}^{\infty}Q_{t}(\cdot,\mathrm{d}u),M\rangle=\int_{s}^{\infty}\frac{\lambda e^{-\lambda(u-t)}\left(\lambda(u-t)\right)^{n-N_{t}}}{(n-N_{t})!}\left(\frac{n-1-N_{t-}}{\lambda(u-t)}-1\right)\lambda\mathrm{d}t,\]
 therefore\[
N_{t}-\int_{0}^{t}\frac{N_{T_{n}-}-N_{u}}{T_{n}-u}\mathrm{d}u,\quad t<T_{n},\]
 is a $\mathbb{G}$-martingale. Note that, by using this proposition,
we cannot extend the $\mathbb{G}$-martingale to values of $t\geq T_{n}.$
\end{example}

\subsection{Progressive enlargement of filtrations}

In the progressive enlargement of filtrations one consider a filtration
$\mathbb{G=}(\mathcal{G}_{t})_{t\in[0,T]}$ with $\mathcal{G}_{t}=\cap_{s>t}(\mathcal{F}_{s}\vee\mathcal{H}_{s})$,
where $\mathbb{H=}(\mathcal{H}_{t})_{t\ge0}$ is another filtration.
The case where $\mathcal{H}_{t}=\sigma(\mathbf{1}_{\{\tau\leq t\}})$
with $\tau$ a random time has been extensely studied, see for instance
(3)-.-, (4) or (1), among others. However, as mentioned in the introduction,
few studies has been developed in the general setting. We present
now an example of an exception, extract from (5), in which $\mathcal{H}_{t}=\sigma(L_{t})$,
$t\in[0,T]$, for $L_{t}=G(X,Y_{t}),$ where $X$ is an $\mathcal{F}_{T}$-measurable
random variable, $(Y_{t})_{t\ge0}$ is a process independent of $\mathcal{F}_{T}$,
and $G$ is a Borelian function\textit{. }\textit{\emph{The following
proposition gives a particular case of this situation, whereas }}in
section {\S}\ref{sub:insider-trading} below, we give the applications
context from which this example arose.

\begin{proposition} \label{Pro: example_peof}Assume that $(B_{t})_{t\in[0,T]}$
is a Brownian motion and take $\mathbb{F}:=\mathbb{F}^{B}$. Let $(W_{t})_{t\in[0,T]}$
be another Brownian motion independent of $(B_{t})_{t\in[0,T]}$,
and consider the process $V_{t}:=B_{T}+\int_{t}^{T}\sigma_{s}\mathrm{d}W_{s}$,
with $\int_{t}^{T}\sigma_{s}^{2}\mathrm{d}s<\infty,$ for all $0\leq t\leq T$.
Then, provided that \[
\int_{0}^{t}\frac{T}{T+\int_{s}^{T}\sigma_{u}^{2}\mathrm{d}u-s}\mathrm{d}s<\infty,\]
 we have that the Doob-Meyer decomposition of $(B_{t})_{t\in[0,T]}$
in $\mathbb{F}^{B,V}$ is given by \[
B_{t}=\tilde{W}_{t}+\int_{0}^{t}\frac{V_{s}-B_{s}}{T+\int_{s}^{T}\sigma_{u}^{2}\mathrm{d}u-s}\mathrm{d}s,\quad0\leq t<T\]
 where $(\tilde{W}_{t})_{t\in[0,T]}$ is a Brownian motion but correlated
with $(V_{t})_{t\in[0,T]}$. \end{proposition}

\begin{proof} $(\tilde{W}_{t})_{t\in[0,T]}$ is a centered Gaussian
processes and for $0\leq s\leq t<T$ \begin{eqnarray*}
\mathbb{E}[\tilde{W}_{t}\tilde{W}_{s}] & = & \mathbb{E}\left[\left(B_{t}-\int_{0}^{t}\frac{V_{u}-B_{u}}{T+\int_{u}^{T}\sigma_{v}^{2}\mathrm{d}v-u}\mathrm{d}u\right)\left(B_{s}-\int_{0}^{s}\frac{V_{u}-B_{u}}{T+\int_{u}^{T}\sigma_{v}^{2}\mathrm{d}v-u}\mathrm{d}u\right)\right]\\
 & = & s-\mathbb{E}\left[B_{t}\int_{0}^{s}\frac{V_{u}-B_{u}}{T+\int_{u}^{T}\sigma_{v}^{2}\mathrm{d}v-u}\mathrm{d}u\right]-\mathbb{E}\left[B_{s}\int_{0}^{t}\frac{V_{u}-B_{u}}{T+\int_{u}^{T}\sigma_{v}^{2}\mathrm{d}v-u}\mathrm{d}u\right]\\
 &  & +\mathbb{E}\left[\int_{0}^{t}\frac{V_{u}-B_{u}}{T+\int_{u}^{T}\sigma_{v}^{2}\mathrm{d}v-u}\mathrm{d}u\int_{0}^{s}\frac{V_{u}-B_{u}}{T+\int_{u}^{T}\sigma_{v}^{2}\mathrm{d}v-u}\mathrm{d}u\right]\\
 & = & s-\int_{0}^{s}\frac{t-u}{T+\int_{u}^{T}\sigma_{v}^{2}\mathrm{d}v-u}\mathrm{d}u-\int_{0}^{s}\frac{s-u}{T+\int_{u}^{T}\sigma_{v}^{2}\mathrm{d}v-u}\mathrm{d}u\\
 &  & +2\mathbb{E}\left[\int_{0}^{s}\left(\int_{0}^{r}\frac{\left(V_{r}-B_{r}\right)\left(V_{u}-B_{u}\right)}{\left(T+\int_{r}^{T}\sigma_{v}^{2}\mathrm{d}v-r\right)\left(T+\int_{u}^{T}\sigma_{v}^{2}\mathrm{d}v-u\right)}\mathrm{d}u\right)\mathrm{d}r\right]\\
 &  & +\mathbb{E}\left[\int_{s}^{t}\left(\int_{0}^{s}\frac{\left(V_{r}-B_{r}\right)\left(V_{u}-B_{u}\right)}{\left(T+\int_{r}^{T}\sigma_{v}^{2}\mathrm{d}v-r\right)\left(T+\int_{u}^{T}\sigma_{v}^{2}\mathrm{d}v-u\right)}\mathrm{d}u\right)\mathrm{d}r\right]\\
 & = & s-\int_{0}^{s}\frac{s+t-2u}{T+\int_{u}^{T}\sigma_{v}^{2}\mathrm{d}v-u}\mathrm{d}u+2\int_{0}^{s}\frac{s-u}{T+\int_{u}^{T}\sigma_{v}^{2}\mathrm{d}v-u}\mathrm{d}u+\int_{0}^{s}\frac{t-s}{T+\int_{u}^{T}\sigma_{v}^{2}\mathrm{d}v-u}\mathrm{d}u\\
 & = & s.\end{eqnarray*}
 On the other hand, for $t\geq s$ \[
\mathbb{E}[\tilde{W}_{t}V_{s}]=s-\int_{0}^{s}\frac{T+\int_{s}^{T}\sigma_{v}^{2}\mathrm{d}v-u}{T+\int_{u}^{T}\sigma_{v}^{2}\mathrm{d}v-u}\mathrm{d}u>0,\]
 provided that $\sigma_{v}$ is not identically null (a.e.). \end{proof}

\begin{remark} It is important to note that contrarily to the case
of initial enlargement, the \textit{innovation process} $(\tilde{W}_{t})_{t\in[0,T]}$
is not necessarily independent of the additional information. Then
this fact makes the application of enlargement of filtrations in our
framework more involved. In other words, in most models, it is assumed
that the privilege information $(V_{t})_{t\in[0,T]}$ is independent
of the demand process of liquidity traders $(\tilde{W}_{t})_{t\in[0,T]}$.
Consequently, the previous proposition cannot be used with these models.
Instead, we have to look for processes such that their Doob-Meyer
decomposition is of the form \[
X_{t}=\tilde{W}_{t}+\int_{0}^{t}\theta(V_{t};X_{u},0\leq u\leq s)\mathrm{d}s,\quad0\leq t\leq T,\]
 where $(\tilde{W}_{t})_{t\in[0,T]}$ and $(V_{t})_{t\in[0,T]}$ are
independent. \end{remark}

Now consider the case when $\mathcal{H}_{t}=$ $\sigma(V_{t})$ for
\[
V_{t}=V_{0}+\int_{0}^{t}\sigma_{s}\mathrm{d}W_{s}^{1},\; t\ge0,\]
 where $\sigma_{s}$ is a deterministic function, $V_{0}$ is a zero
mean normal r.v., and $(W_{t}^{1},W_{t}^{2})_{t\in[0,T]}$ is a 2-dimensional
Brownian motion independent of $V_{0}$. We have the following proposition:

\begin{proposition} \label{prog}Assume that $Var(V_{T})=1$and that
\[
\int_{0}^{t}\frac{\mathrm{d}s}{Var(V_{s})-s}<\infty\text{ for all }0\leq t<T.\]
 Then \[
B_{t}=W_{t}^{2}+\int_{0}^{t}\frac{V_{s}-B_{s}}{Var(V_{s})-s}\mathrm{d}s,\;0\leq t\leq T,\]
 is a Brownian motion with $B_{T}=V_{T}.$ \end{proposition}

\begin{proof} Let us denote $v_{x}:=Var(V_{x})$, for $0\leq x\leq T$.
We have \[
B_{t}=\int_{0}^{t}\exp\left(-\int_{y}^{t}\frac{1}{v_{x}-x}\mathrm{d}x\right)\mathrm{d}W_{y}^{2}+\int_{0}^{t}\exp\left(-\int_{u}^{t}\frac{1}{v_{x}-x}\mathrm{d}x\right)\frac{V_{y}}{v_{y}-y}\mathrm{d}y,\]
 so $(B_{t})_{t\in[0,T]}$ is a centered Gaussian process, and for
$s\leq t<T$,\begin{eqnarray*}
\mathbb{E}[B_{t}B_{s}] & = & \exp\left(-\int_{s}^{t}\frac{1}{v_{x}-x}\mathrm{d}x\right)\\
 &  & +\mathbb{E}\left[\int_{0}^{t}\int_{0}^{s}\exp\left(-\int_{y}^{t}\frac{1}{v_{x}-x}\mathrm{d}x\right)\exp\left(-\int_{z}^{s}\frac{1}{v_{x}-x}\mathrm{d}x\right)\frac{V_{y}V_{z}}{(v_{y}-y)(v_{z}-z)}\mathrm{d}y\mathrm{d}z\right]\\
 & = & \exp\left(-\int_{s}^{t}\frac{1}{v_{x}-x}\mathrm{d}x\right)\int_{0}^{s}\exp\left(-2\int_{y}^{s}\frac{1}{v_{x}-x}\mathrm{d}x\right)\mathrm{d}y\\
 &  & +\int_{s}^{t}\int_{0}^{s}\exp\left(-\int_{y}^{t}\frac{1}{v_{x}-x}\mathrm{d}x\right)\exp\left(-\int_{z}^{s}\frac{1}{v_{x}-x}\mathrm{d}x\right)\frac{v_{z}}{(v_{y}-y)(v_{z}-z)}\mathrm{d}y\mathrm{d}z\\
 &  & +2\int_{0}^{s}\int_{0}^{y}\exp\left(-\int_{u}^{t}\frac{1}{v_{x}-x}\mathrm{d}x\right)\exp\left(-\int_{z}^{s}\frac{1}{v_{x}-x}\mathrm{d}x\right)\frac{v_{z}}{(v_{y}-y)(v_{z}-z)}\mathrm{d}y.\end{eqnarray*}
 Then, since \[
\int_{0}^{s}\exp\left(-\int_{z}^{s}\frac{1}{v_{x}-x}\mathrm{d}x\right)\frac{v_{z}}{v_{z}-z}\mathrm{d}z=s,\]
 and \[
2\int_{0}^{s}\exp\left(-2\int_{z}^{s}\frac{1}{v_{x}-x}\mathrm{d}x\right)\frac{v_{z}}{v_{z}-z}\mathrm{d}z=2s+\int_{0}^{s}\exp\left(-2\int_{y}^{s}\frac{1}{v_{x}-x}\mathrm{d}x\right)\mathrm{d}y\]
 we obtain that $\mathbb{E}[B_{t}B_{s}]=s.$ So for $0\leq t<T$ we
have that $(B_{t})_{t\in[0,T]}$ is a standard Brownian motion. On
the other hand \begin{eqnarray*}
\mathbb{E}[B_{t}V_{t}] & = & \mathbb{E}\left[\int_{0}^{t}\exp\left(-\int_{y}^{t}\frac{1}{v_{x}-x}\mathrm{d}x\right)\frac{V_{y}V_{t}}{v_{y}-y}\mathrm{d}y\right]\\
 & = & \int_{0}^{t}\exp\left(-\int_{y}^{t}\frac{1}{v_{x}-x}\mathrm{d}x\right)\frac{v_{y}}{v_{y}-y}\mathrm{d}y\\
 & = & t,\end{eqnarray*}
 therefore\begin{eqnarray*}
\mathbb{E}[(B_{t}-V_{t})^{2}] & = & \mathbb{E}[B_{t}^{2}]+\mathbb{E}[V_{t}{}^{2}]-2\mathbb{E}[B_{t}V_{t}]\\
 & = & t+v_{t}-2t=v_{t}-t,\end{eqnarray*}
 and, since by hypothesis $v_{T}=1,$ this means that \[
\lim_{t\rightarrow T}B_{t}\overset{L^{2}}{=}V_{T},\]
 then for all $0\leq t<T$ \[
\mathbb{E}\left[\int_{0}^{t}\frac{|V_{s}-B_{s}|}{v_{s}-s}\mathrm{d}s\right]<\int_{0}^{t}\frac{\mathbb{E}[(V_{s}-B_{s})^{2}]^{\frac{1}{2}}}{v_{s}-s}\mathrm{d}s=\int_{0}^{t}\sqrt{v_{s}-s}\mathrm{d}s<\sqrt{2},\]
 and this implies, by the monotone convergence theorem, that \[
\lim_{t\rightarrow T}\int_{0}^{t}\frac{|V_{s}-B_{s}|}{v_{s}-s}\mathrm{d}s=\int_{0}^{T}\frac{|V_{s}-B_{s}|}{v_{s}-s}\mathrm{d}s<\infty\text{ }\]
 and that $B_{T}=\lim_{t\rightarrow T}B_{t}$ is well defined. Now,
we have, by the uniqueness of the limit in probability, that $V_{T}=B_{T}$
a.s. \end{proof}

\section{Applications of Enlargement of Filtrations to Mathematical Finance}

In mathematical finance, the enlargement of filtration theory may
be applied to credit risk and to insider trading. In credit risk,
the original filtration may be thought to represents the information
related to defaultable-free assets in the market. We then construct
the enlarged filtration by introducing the information related to
a defaultable prone asset of interest. On the other hand, for inside
trading, the reference filtration is that of a regular market agent.
This filtration is enlarged by the \emph{insider} agent who has privileged
information about, for instance, the future value of said asset. In
the following subsections we give more details about said applications.

\subsection{Applications to credit risk theory\label{sub:credit-risk}}

The main objective of quantitative models of credit risk is to provide
ways to price and hedge financial contracts that are sensitive to
\emph{credit risk}, that is to say, the risk of an economic loss due
to the failure -or \emph{default}- of a counterpart to fulfill its
contractual obligations. When a contract or a firm defaults it is
said that the \emph{default event} occurs, and the random time $\tau$
at which the default event occurs is called \emph{default time}. A
vast majority of mathematical research devoted to credit risk is concerned
with modelling default times. Two methodologies have emerged in order
to model the default time: the \emph{structural approach} which dates
back to Black and Scholes (cf. (6)) and Merton (cf. (7)), and the
\emph{reduced-form approach} originated with Jarrow and Turnbull (cf.
(8)).

Within structural models the value of the firm of interest is assumed
to have the following dynamics under the neutral probability $\mathbb{P}^{*}$\[
V_{t}:=\exp(L_{t}),\quad t\in[0,T],\]
 where $(L_{t})_{t\in[0,T]}$ is a Lévy process. More over, credit
events are triggered by movements of the firm's value $(V_{t})_{t\in[0,T]}$
to some random or not-random lower threshold, sometimes called \emph{default
barrier}. For instance, taking a constant barrier $K$, the default
time is defined as\[
\tau:=\inf\{t\in[0,T]\;:\; V_{t}<K\},\]
 where $K<V_{0}$, and $\tau$ is set to be $\infty$ if $(V_{t})_{t\in[0,T]}$
does not crosses the default barrier. Notice that the default event
is defined endogenously within the model, and the information available
to the modeler has to be the same that the firm's manager has. More
details about this approach can be found references such as (9), and
(10).

In the reduced-form approach, the modeler does not have the full information
that the firm manager possesses but only a subset of it, generated
by the \emph{default process} $(H_{t}:=\mathbf{1}_{\{\tau<t\}},\; t\in[0,T])$,
and several other related state variables. The value of the firm's
assets and its capital is not modelled at all, and credit events are
specified in terms of some exogenously specified jump process. In
the literature, reduced form framework has been split into two different
approaches, the \emph{Hazard Process Models} and the \emph{Intensity-Based
Models}, depending on whether the information of the default free
assets is introduced or not.

In Hazard Process models, a filtration $\mathbb{F}$ (generated by
a Brownian motion) is interpreted as the information related to the
default-free assets of the market. Let us denote by $\mathbb{H}:=(\mathcal{H}_{t})_{t\in[0,T]}$
the filtration generated by the default process $(H_{t}:=\mathbf{1}_{\{\tau<t\}},\; t\in[0,T])$.
We can consider a filtration $\mathbb{G}:=(\mathcal{G}_{t})_{t\in[0,T]}$
encompassing the information regarding to default-free assets, as
well as the information regarding the defaultable asset of interest.
More specifically, we consider the progressive enlargement given by
$\mathcal{G}_{t}:=\mathcal{F}_{t}\vee\mathcal{H}_{t}$, for every
$t\in[0,T]$. It should be emphasized that $\tau$ is not necessary
a stopping time with respect to $\mathbb{F}$, though it is a stopping
time with respect to $\mathbb{G}$. It is well known that in order
to preclude arbitrage opportunities in a default-free market, the
properly discounted asset prices have to be $\mathbb{F}$-semimartingales.
Since the full market is also assumed to be arbitrage free, these
prices must be $\mathbb{G}$-semimartingales as well.

Although the structural and the reduced-from approaches for credit
risk modelling seem to be conceptually different, efforts has been
made in order to establish relationships between them. We distinguish
between two main lines of work in this matter: (\emph{i}) creation
of a general model for credit risk modelling encompassing the two
approaches; and (\emph{ii}) the pass from one approach to the other
by modifying the information available to the modeler. See for instance,
(11) and (12) for approach (\emph{i}), and (13) for approach (\emph{ii}).

A different way to relate the structural and reduced-form models may
be done by considering a Kyle-Back model for insider trading (cf.
(14) and (15)). Following this approach, (16) presents a model of
\emph{asymmetric information} (\emph{i.e.}, when different market
agents possess different informations about the market) in which both
approaches play a role.

\subsection{Applications to insider trading\label{sub:insider-trading} in a
Kyle-Back market model}

A company issues a risky asset. Assume that the process $(Z_{t})_{t\in[0,T]}$
models the value of the company. This process may be taken as a Brownian
motion, or as a more general process that may have a drift and jumps.
Three types of agents interact in the Kyle-Back market we are considering:
\begin{itemize}
\item The \emph{noise traders} who trade for liquidity or hedging reasons.
They observe only their own cumulative demands -modelled by Brownian
motion $(B_{t})_{t\in[0,T]}$, started at $0$, and independent from
$(Z_{t})_{t\in[0,T]}$- and whether the risky asset has defaulted
or not.
\item The \emph{informed trader} who is an agent that observes continuously
in time the defaultable bond prices, and knows some additional information
about the risky asset (\emph{e}.\emph{g}., its price at some prefixed
time). Let the random variable $X$ contains the privileged. It is
plausible to think that the insider does not knows exactly $X$ but
a good estimation of it, say for instance, that she knows $X$ modified
by some perturbing noise. More specifically, let the additional information
until time $t$ be given by a family of random variables $(L_{s})_{s\leq t}$.
We suppose that these random variables have the following structure\[
L_{t}=G(X,Y_{t}),\]
 where $X$ is an $\mathcal{F}_{T}$-measurable random variable, the
process $(Y_{t})_{t\in\left[0,T\right]}$ is independent of the $\sigma$-algebra
$\mathcal{F}_{T}$, and $G:\mathbb{R}^{2}\rightarrow\mathbb{R}$ is
a given measurable function. In that sense we can define the insider's
filtration $\mathbb{F}_{I}$ as the filtration $\mathbb{F}^{B}$ enlargement
by $\mathbb{F}^{L}$, that is to say, $\mathbb{F}_{I}:=(\cap_{s>t}(\mathcal{F}_{s}^{B}\vee\mathcal{F}_{s}^{L}),\; t\in\left[0,T\right])$.
The random variables $Y_{t}$ represent the additional noise, whereas
the function $G$ specifies how the perturbation is made. One expects
in general that $Y_{T}=0$ and that the variance of the noise should
decrease to zero as time approaches the moment at which the additional
information (possed by the insider) is released to the public.
\item The \emph{market maker} is that one market agent that observes\emph{\
total order} $(R_{t})_{t\in[0,T]}$ of the noise traders and the insider,
and sets the price of the risky asset. Let $\mathbb{F}^{R}$ denotes
the minimal right continuous and complete filtration generated by
$(R_{t})_{t\in[0,T]}$. Consequently, the market maker's information
$\mathbb{F}^{M}$ is given by the progressive filtration enlargement
$(\cap_{s>t}(\mathcal{F}_{s}^{Y}\vee\sigma(\tau\wedge s)),\; t\in[0,T])$. 
\end{itemize}
A particular case of this scenario has been seen with Proposition
\ref{Pro: example_peof} where it is considered\[
L_{t}:=B_{T}+\int_{t}^{T}\sigma_{s}\text{d}W_{s},\]
 being $(W_{t})_{t\in\left[0,T\right]}$ a Brownian motion, independent
of $(B_{t})_{t\in\left[0,T\right]}$, and $\int_{t}^{T}\sigma_{s}\text{d}s<\infty$.

In (6), the authors  study the case of a company issuing a defaultable
bond, with face value $1$, and maturity time $T=1$. For simplicity,
the interest rate is taken as zero, and the value of the company is
assumed to follow a Brownian motion $(Z_{t})_{t\in[0,T]}$ . In turn,
default is set as\[
\tau=\inf\{t\in[0,T]\;:\; Z_{t}=-1\}.\]
 In this scenario, the privileged information is the default time
$\tau$ and no perturbation is considered, so that the insider's information
is given by the initial filtration enlargement $\mathbb{F}_{I}=(\cap_{s>t}(\mathcal{F}_{s}^{B}\vee\sigma(\tau)),\; t\in[0,T])$.

After properly defining the set of insider's trading strategies $\mathcal{A}$,
and the market maker's pricing rules $\mathcal{H}$, the idea is to
find pairs $(H,\theta)\in\mathcal{H}\times\mathcal{A}$ such that
both the insider and the market market fulfill their respective objectives.
Such a pair is called an \emph{equilibrium. }The insider's objective
is to maximize her expected wealth at time $T$, provided that she
is risk-neutral. Whereas the market maker's objective is to set a
\emph{rational} price of the risky asset in order to \emph{clear the
market}. The main result in (16) is the existence of an equilibrium
$(H^{*},\theta^{*})$ for which the $\tau$ is a predictable stopping
time under the market maker's information $\mathbb{F}^{M}$, and such
that the equilibrium total order solves $(R_{t}^{*})_{t\in[0,T]}$
the following stochastic differential equation (SDE)\[
\mathrm{d}R_{t}=\mathrm{d}B_{t}+\left(\frac{1}{1+R_{t}}-\frac{1-R_{t}}{\tau-t}\right)\mathbf{1}_{\{\tau\leq t\}}\mathrm{d}t.\]
 Aside from an equilibrium characterization lemma proved therein,
the proof of this result relays on the construction of a weak solution
to the SDE. As seen in Example \ref{G4}, $(Z_{t})_{t\in[0,T]}$ has
the following $\mathbb{F}_{I}$-decomposition. \[
\mathrm{d}Z_{t}=\mathrm{d}\beta_{t}+\left(\frac{1}{1+Z_{t}}-\frac{1-Z_{t}}{\tau-t}\right)\mathbf{1}_{\{\tau\leq t\}}\mathrm{d}t,\]
 where $(\beta_{t})_{t\in[0,T]}$ is a $\mathbb{F}_{I}$-Brownian
motion. It can be proved that the SDE possesses a unique strong solution,
and two consequences follow from this. On the one hand, $(R_{t}^{*})_{t\in[0,T]}$
has the same law as $(Z_{t})_{t\in[0,T]}$, and thus $(R_{t}^{*})_{t\in[0,T]}$
is a Brownian motion in its own filtration. On the other hand, $\tau=\inf\{t\in[0,T]\;:\; R_{t}^{*}=-1\}$.
Hence $\tau$ is a stopping time with respect to $\mathbb{F}^{R^{*}}$,
and the filtrations $\mathbb{F}^{R^{*}}$ and $\mathbb{F}^{M}$coincide,
so that $(R_{t}^{*})_{t\in[0,T]}$ is an $\mathbb{F}^{M}$-Brownian
motion, too.

\end{document}